Sandra Bella
Postdoc ERC Philiumm adg n°101020985, Laboratoire Sphère, Université de Paris-Cité
bellusky@hotmail.com

# The Peculiar Destiny of *Sentiment de Monsieur Leibnitz* (May 1705 – March 1706)[1]

## 1. Introduction

In July 1700, at the Académie royal des sciences, Michel Rolle (1652-1719) launched his criticism about the soundness of Leibniz' new differential calculus. The defence was assured by Pierre Varignon (1654-1722). This episode was the beginning of the well-known "querelle des infiniment petits", which did not calm until 1706. However, this quarrel was contrary to the academic rules, which required that its members should exchange "avec ménagement[2]." The Academy requested silence[3]. This implied that any polemic memoir concerning differential calculus was forbidden to be lectured during the academic sessions. The real place of the debate was shifted. The *Journal des Sçavans* now published articles by Rolle and the mathematician Joseph Saurin (1659-1737) who publicly relayed Varignon as the herald of the defence of differential calculus.

Leibniz was deeply concerned. Throughout the quarrel, he composed writings to publicly justify his calculus to the Parisian scholars[4]. All these texts have been the object of historical and philosophical analyses[5]. But only three of these texts were actually published, and of these, *Sentiment de Monsieur Leibnitz* was, as will be shown, immediately confiscated.

---

[1] This article is a shortened version of the article « *Magis morale quam mathematicum*, L'attestation volée (mai 1705 - mars 1706) », *Studia Leibnitiana* 51, 2019/2, 176–202.

[2] Among the fifty articles of the new 1699 regulation, the twenty-sixth is dedicated to a rule of conduct between members, especially when they are of differing opinions. "L'Académie veillera exactement a ce que dans les occasions ou quelques académiciens seront d'opinions différentes, ils n'employent aucuns termes de mépris ni d'aigreur l'un contre l'autre, soit dans leurs discours, soit dans leurs écrits ; et lors même qu'ils combatteront les sentimens de quelques Sçavans que ce puisse être, l'Académie les exhortera à n'en parler qu'avec ménagement.", *Registres des Procès-verbaux de l'Académie Royale des sciences* (from now on "*PVARS*", followed by the tome, and the folio), 18, 110 v°.

[3] "comme la dispute des infiniment petits traînait trop en longueur Mr l'Abbé Bignon a nommé pour commissaires devant qui tout se passera, le P.Gouye, MMrs Cassini et de la Hire", *PVARS,* 20, 335v°.

[4] Gottfried Wilhelm Leibniz: "Mémoire de Mr Leibnitz touchant son sentiment sur le Calcul différentiel », in : *Journal de Trévoux*, novembre – décembre 1701, pp. 270–272, GM IV, 95–96, GM V, 350 ; "Extrait d'une lettre de M. Leibnitz à M. Varignon, contenant l'explication de ce qu'on a rapporté de luy dans les Memoires des mois de Novembre et Decembre derniers", in: *Journal des Sçavans*, 20 mars 1702, pp. 183–186, GM IV, 91–95 ; "Justification du Calcul des infinitésimales par celuy de l'Algèbre ordinaire", LBr 951, f°15–16, GM IV, 104–106 ; « Cum prodiisset », LH XXXV, 7, 19, fols. 1–6, Carl Immanuel Gerhardt : *Historia et origo calculi differentialis a G. G. Leibnitio conscripta*, Hannover 1846, pp. 39–50 ; « Defense du calcul des Differences », LH XXXV, 7, 22, f°1–2, Appendice de Enrico Pasini : *La nozione di infinitesimo in Leibniz : Tra matematica e metafisica*, diss., Torino, Università degli studi, Fac. di Lettere e filosofía, 1985–1986, pp. 14–19 ; "Quaestio de jure negligendi quantitates infiniti parvas", LH XXXV, 7, 21, f° 1–2, Appendice de Pasini : *La nozione* (cf. note 4), pp. 40–47 ; "Puisque des personnes que j'estime", LH XXXV, 7, 9, f° 1–2 ; "Sentiment de Monsieur Leibnitz", Gotha FB A 448–449, f° 41–42 ; Leibniz: "Sentiment de Monsieur Leibnitz", in: Joseph Saurin : *Continuation de la défense de M. Saurin contre la Réplique de M. Rolle publiée en 1703, sous le titre de Remarques touchant le problème général des Tangentes, etc.*, Amsterdam 1706, pp. 35–39. Some of these manuscripts are presented and transcribed at the following web page: https://eman-archives.org/philiumm/rubric/justification-of-the-differential-calculus (see on 16 May 2023).

[5] Richard T. W. Arthur: "Leery Bedfellows : Newton and Leibniz on the Status of Infinitesimals" in: Ursula Goldenbaum / Douglas Jesseph (eds.): *Infinitesimal Differences. Controversies between Leibniz and His Contemporaries*, Berlin 2008, pp. 7–30 ; Philip Beeley: "Infinity, infinitesimals, and the Reform of Cavalieri: John Wallis and his critics", in: Goldenbaum/Jesseph: *Infinitesimal Differences* (cf. ci-dessus), pp. 31–52 ; Henk J. Bos: "Differentials, Higher-Order Differentials and the Derivative in the Leibnizian Calculus", in: *Archive for History of Exact Sciences* 14 (1974), pp. 2–90 ; Douglas Jesseph: "Leibniz on the elimination of infinitesimals", in: Norma B. Goethe / Philip Beeley / David Rabouin (eds.):

"Sentiment de Monsieur Leibnitz" and its short version in French "Puisque des personnes que j'estime" are two texts to be situated at the peak of the *querelle des infiniment petits* in 1705. In his essay *Remarques de M. Rolle de l'Académie Royale des Sciences touchant le problesme general des tangentes*[6], Rolle undermined the foundations of differential calculus and in particular showed the insufficiency of calculus to deal with the problem of determining the tangent of a quartic curve. Two years later, on 23 April 1705, Saurin finally published a reply[7]. In this article, Saurin implored the Academy to judge the discordances between himself and Rolle, which he listed in seven points[8]. This extremely tense climate prompted Varignon to write to Leibniz on 10 May 1705[9]. He hoped that Leibniz would intervene with the Academy. He also asked him to solicit the members of the République des lettres who *entendent* his calculation to obtain "attestations" confirming the exactitude of the contents treated by Saurin in his article. This episode has been the subject of studies, albeit without access to these two manuscripts, nor to the corresponding epistolary exchanges.

Leibniz promptly reacted. He wrote two letters, one to Jean Gallois (1632-1707), who, together with Rolle, was one of the most opposed personalities to his calculus in the Academy, and the other to Jean-Paul Bignon (1662-1743), President of the Academy. He therefore wrote the requested attestation. Two versions exist: a holographic Latin version beginning with the words *"Cum viri docti et a me aestimati"* (this manuscript is conserved at the Forschungsbibliothek Gotha, from now on Gotha FB), and an abridged version in French, in the hand of a copyist, registered under the title "Puisque des personnes que j'estime". Only the Latin version was finally published under the title "Sentiment de Monsieur Leibnitz". It is joined to a text by Saurin[10] – which is a continuation of his article of 23 April 1705 – and to other attestations: that of Jacob Hermann (1678-1733) and those of the Bernoulli brothers. However, all the publications were confiscated by Bignon and perhaps partly destroyed. It was only recently that I identified one of the copies in the British Library[11]. It is catalogued by the title corresponding to Saurin's writing but unfortunately without any indication of the authors of the attestations and in particular of Leibniz. It is probably for this reason that no copy was considered to have survived the Bignonian authority.

As mentioned above, "Sentiment de Monsieur Leibnitz" is one of a group of texts that Leibniz wrote in order to justify his calculus to the Parisian scholars. For a long time

---

*G. W. Leibniz. Interrelations Between Mathematics and Philosophy*, Netherlands 2015, pp. 189–206 ; Mikhail Katz / David Sherry: "Leibniz's Infinitesimals: Their Fictionality, Their Modern Implementations, and Their Foes from Berkeley to Russell and Beyond", in: *Erkenntnis* 78/3 (2013), pp. 571–625 ; Pasini: *La nozione* (cf. note 4) ; Enrico Pasini : "Die private Kontroverse des G. W. Leibniz mit sich selbst. Handschriften über die Infinitesimalrechnung im Jahre 1702", in : *Leibniz : Tradition und Aktualität. V. Internationaler Leibniz-Kongreß, Hannover, 14.–19. November 1988*, 2 vols., Hannover 1988–1989, vol. 1: *Vorträge. Teil 1*, Hannover 1988, pp. 695–709 ; Enrico Pasini : *Il reale e l'immaginario. La fondazione del calcolo infinitesimale nel pensiero di Leibniz*, Torino 1993 ; Paolo Mancosu: *Philosophy of Mathematics and Mathematical Practice in the Seventeenth Century*, New York 1996 ; David Rabouin / Richard T. W. Arthur: "Leibniz's Syncategorematic Infinitesimals II: Their Existence, Their Use and Their Role in the Justification of Differential Calculus", in: *Archive for History of Exact Sciences* 74 (2020), pp. 401–443.
[6] Michel Rolle: *Remarques de M. Rolle de l'Académie des Sciences touchant le problesme général des tangentes*, Paris 1703.
[7] Joseph Saurin: "Défense de la réponse à M. Rolle de l'Académie royale des sciences contenue dans le Journal des Sçavans du 3 août 1702 CONTRE La réplique de cet auteur publiée en 1703 sous le titre de Remarques touchant le Problème général des Tangentes", in: *Journal des Sçavans*, jeudi 23 avril 1705, pp. 241–256.
[8] *Ibid.*, p. 254.
[9] Letter from Pierre Varignon to Leibniz, 10 May [1705], A III, 9, 549–550 (Vorausedition, 2.11.2018).
[10] Cf. footnote 3.
[11] Under the reference: BLL01017362909.

considered lost, the publication of this manuscript is now an important element for a complete understanding of the *querelle* and its outcome[12].

In the light of this unpublished manuscript, but also of recently published letters addressed to prestigious members of the Academy, our contribution aims to reconstitute the mathematical and political stakes within the Académie royale des sciences, of one of the latest episodes of the *querelle des infiniment petits*.

2. Leibniz's Calculus in the *Académie royale des sciences*

Leibniz's calculus was introduced in Paris as early as December 1691 by Johann Bernoulli (1667-1748) and a circle of mathematicians around the philosopher Nicolas Malebranche (1638-1715)[13]. Among them, Guillaume de L'Hospital (1661-1704) was the first to benefit from Johann Bernoulli's services – private lessons and a rich epistolary correspondence. As a result, he was the first to adopt the new calculus, which he publicised first by reading memoirs at the Académie royale, a member of which he had been since 1693, and then by publishing the first treatise on differential calculus, entitled *Analyse des infiniment petits pour l'intelligence des lignes courbes*[14] in 1696. L'Hospital's treatise thus became the main reference on the subject for the Académie royale.

In 1699, new rules were established at the Academy under the impulse of its President Bignon[15]. These rules had an impact on the diffusion of differential calculus. Many mathematicians from Malebranche's circle were elected. Among the members of the *République des lettres*, who were versed in the new calculus, Leibniz and the Bernoulli brothers were elected as foreign members. Also, article twenty-one stated that at the beginning of the year, every academician (except the honorary members) had to announce the purpose of their research[16]. Varignon, Gallois and Rolle were all pensioners[17]. In return for their pension granted by the King, they had to report on their work to their peers. In the camp of those opposed to differential calculus, Gallois announced that he would revive the geometry of the Ancients because he considered that it had been neglected in favour of "the new methods[18]". By "new methods", he was referring principally to the differential calculus. Finally, Rolle planned to use and develop the algebra to provide a "théorie réglée des méthodes que l'algèbre ait fournies" in order to apply it to geometry.

---

[12] For analyses of the dispute, see, in particular Mancosu: *Philosophy of Mathematics* (cf. footnote 4), pp. 173–178 ; et John B. Shank: *Before Voltaire. The Origin of "Newtonian" Mechanics, 1680–1715*, Chicago/London 2018, pp. 286–297.
[13] Not all are Parisians and some reside in Normandy (Rouen, Angers, Honfleur). Louis Byzance (1647-1722), Louis Carré (1647-1711), François de Catelan (floruit 1676-1710), Bernard Le Bouyer de Fontenelle (1657-1757), Claude Jaquemet (1651-1729), Bernard Lamy (1640-1715), Thomas Fantet de Lagny (1660-1734), Jacques Lelong (1665-1721), Jean Prestet (1648-1690), Charles Reyneau (1656-1728), Pierre Varignon (1654-1722).
[14] Guillaume de L'Hospital: *Analyse des infiniment petits, pour l'intelligence des lignes courbes*, Paris 1696
[15] For a discussion of the regulations, see Shank : *Before Voltaire* (cf. footnote 11), pp. 217–224 ; and David J. Sturdy: *Science and Social Status. The Members of the Académie des Sciences, 1666–1750*, Woodbridge 1995, pp. 281–350.
[16] "Chaque académicien pensionnaire sera obligé de déclarer par écrit à la Compagnie le principal ouvrage auquel il se proposera de travailler : et les autres académiciens seront invités à donner une semblable déclaration de leurs desseins.", *PVARS* 18, 109rº.
[17] Rolle is an algebrist, the others are nominated as geometers.
[18] "Comme presque tous les géomètres ne s'appliquent maintenant qu'aux nouvelles méthodes de géométrie ce que Mr l'Abbé Gallois juge qu'il est important de ne pas abandonner entièrement la manière des Anciens, il prend pour son partage d'illustrer les ouvrages des Anciens géomètres", *PVARS*, 23, 133rº.

This variegated composition did not obviously create a climate favourable to dialogue, but rather revealed tensions that conduced to the eruption of the *querelle des infiniment petits* in July 1700[19].

3. Around a double point

During the quarrel, Rolle submitted several examples of curves for which he considered that the application of the Leibnizian calculus did not enable to find their tangents in a satisfactory way. Among these examples, the equation of quartic curve

$$y^4 - 8y^3 - 16y^2 - 12xy^2 + 48xy + 4x^2 - 64x = 0 \quad (1)$$

was the focus of mathematicians' discussions as it presents a double point (with coordinates x=2 and y=2) for which exist two tangents.

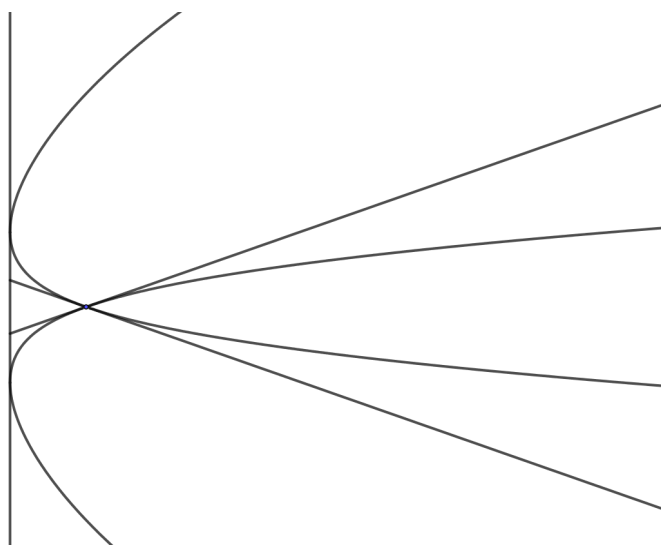

*Figure 1 :* Quartic $y^4 - 8y^3 - 16y^2 - 12xy^2 + 48xy + 4x^2 - 64x = 0$ et the two tangents in the double point
(figure done by the Author)

On April 13, 1702, Rolle published a memoir "Règles et Remarques, pour le problème général des tangentes[20]" in which he explained that the "ordinary methods" – without specifying which ones – were not appropriate for determining any of the subtangents at the double point[21]. He provided without justification an algebraic algorithm for determining the values of the two subtangents being sought[22]. Understanding this allusion as an offensive against the differential method, Varignon, L'Hospital and Saurin wanted to take up this challenge and thus illustrate the fruitfulness of the Leibnizian algorithm.

---

[19] For studies on the diffusion of differential calculus in the Academy of Sciences, see: Sandra Bella: *La (Re)construction de l'analyse infinitésimale de Leibniz (1692-1706)*, Garnier, 2022; Michel Blay: *La naissance de la mécanique analytique. La science du mouvement au tournant des XVII<sup>e</sup> et XVIII<sup>e</sup> siècles*, Paris 1992 ; Mancosu: *Philosophy of Mathematics* (cf. footnote 4) ; and Shank : *Before Voltaire* (cf. footnote 11).

[20] Michel Rolle : « Règles et remarques, pour le problème général des tangentes », in: *Journal des Sçavans*, 13 avril 1702, pp. 239–253.

[21] The subtangent at a point $M$ is the algebraic length between the intersection of the tangent and the $x$-axis and the orthogonal project of $M$ on the $x$-axis (cf. *Figure 1*).

[22] The idea of the algorithm is to replace $y$ by $y + nz$ and $x$ by $x + nv$ in the quartic equation. The algorithm computes the terms of the new equality "in slices". The first slice (which corresponds to the determination of the tangents at a non-singular point) is formed from the terms of the equation of the curve by multiplying each of the $y$ terms by their exponent and then replacing one of the $y$ terms in the product with $z$. For example $y^4$ is transformed into $4y^3z$. The same type of procedure is applied to the $x$-terms but replacing by $v$. Then all these products are multiplied by $n$. For the second block, a similar procedure is followed from the first block, obtained in the first step, except that this time the products are multiplied by $n^2$. Thus, the algorithm calculates a block from the previous one in retrograde order. Rolle presents the new equality arranged in order. In the case of the quartic, by replacing $x$ and $y$ by 2, the first part of the equality thus arranged annuls, so we must move on to the next block which this time does not annul and provides the values of the two subtangents, Rolle: *Remarques* (cf. footnote 5), pp. 240-242.

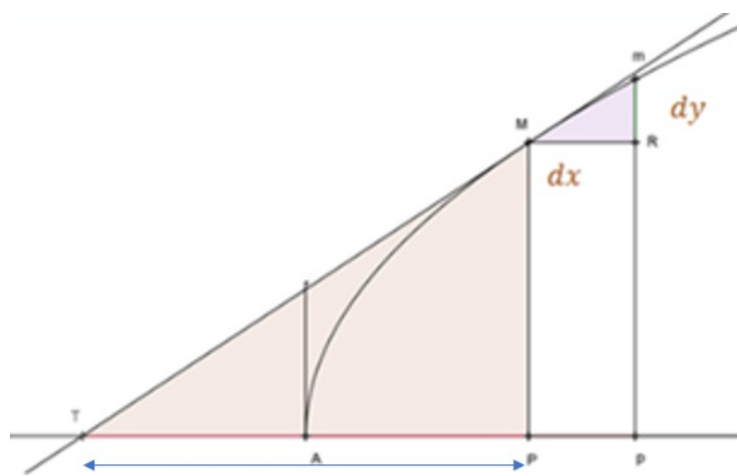

*Figure 2: Characteristic triangle of sides $dx$ and $dy$ similar to $TMP$ (Figure done by the Author)*

L'*Analyse des infiniment petits pour l'intelligence des lignes courbes*[23] provides a general method for determining the subtangent $t$ at a point $M$ of a curve by providing the formula $t = y\frac{dy}{dx}$, which is based on geometrical considerations (*Figure 2*)[24]. But, in the case of the double point, this formula is not directly applicable.

Indeed, by differentiating the equality
$$y^4 - 8y^3 - 16y^2 - 12xy^2 + 48xy + 4x^2 - 64x = 0 \quad (1)$$
then
$$\frac{dy}{dx} = \frac{3y^2 - 12y - 2x + 16}{y^3 - 6y^2 + 8y - 6xy + 12x} \quad (2)$$
and with $x = 2$ and $y = 2$, the numerator and denominator become zero. Article 163 provides a solution to this impasse: one should differentiate the numerator and denominator before replacing them with the values that previously annulled them[25].

Saurin, assisted by Varignon and probably other members of the circle around Malebranche[26], replied to Rolle in an article published in the *Journal des Sçavans* on 3 August 1702[27]. When $x$ and $y$ are equal to 2 $\frac{dy}{dx}$ is indeterminate. To determine this value, it is appropriate, according to article 163, to differentiate the numerator and the denominator of
$$\frac{dy}{dx} = \frac{3y^2 - 12y - 2x + 16}{y^3 - 6yy + 8y - 6xy + 12x} \quad (2)$$
then
$$\frac{dy}{dx} = \frac{d(3y^2 - 12y - 2x + 16)}{d(y^3 - 6yy + 8y - 6xy + 12x)} = \frac{6ydy - 12dy - 2dx}{3yydy - 12ydy + 8dy - 6xdy - 6ydx + 12dx} \quad (3)$$

---

[23] L'Hospital: *Analyse des infiniment petits* (cf. footnote 13), p. 11.
[24] See *Figure 2*: Around of point $M$, considering a point $m$ infinitely close to $M$, the arc $Mm$ is assumed to be rectilinear. The curvilinear triangle $MmR$ is thus assumed to be rectilinear and similar to $TMP$, hence the equality.
[25] "Soit une ligne courbe $AMB$ ($AP = x, PM = y, AB = a$) telle que la valeur de l'appliquée $y$ soit exprimée par une fraction, dont le numérateur et dénominateur deviennent chacun zero lorsque $x = a$, c'est-à-dire lorsque le point $P$ tombe sur le point $B$. On demande quelle doit être la valeur de l'appliquée $BD$. (…) si l'on prend la différence du numérateur, & qu'on la divise par la différence du dénominateur, après avoir fait $x = a = Ab$ ou $AB$, l'on aura la valeur de l'appliquée $BD$ ou $bd$.", L'Hospital: *Analyse des infiniment petits* (cf. note 13), p. 145.
[26] Bella: *La (Re)construction* (cf. footnote 18), pp. 436-444.
[27] Joseph Saurin: "Reponse à l'ecrit de M. Rolle de l'Ac. R. des Sc. inseré dans le Journal du 13 avril 1702. sous le titre de Regles et Remarques pour le Problème general des Tangentes", in: *Journal des Sçavans*, jeudi 3 août 1702, pp. 519–534.

Substituting $x$ by 2 and $y$ par 2, the equality (3) becomes $\frac{dx}{8dy} = \frac{dy}{dx}$ (4)

By multiplying the extremes:
$$\frac{dy^2}{dx^2} = \frac{1}{8}. \text{ (5)}$$

As $x = 2$, (5) becomes $t^2 = \frac{x^2 dy^2}{dx^2} = \frac{4}{8}$. (6) then $t = \frac{xdy}{dx} = \pm\frac{\sqrt{2}}{2}$.

Rolle did not approve this way of proceeding despite this being the same result as that he obtained with his algebraic method. First of all, he explained that "on a fait concourir deux Methodes de cette Analyse, l'une de l'art. 9 l'autre de l'art. 163[28]", so Saurin's way was not a general method for determining the tangents of curves. Moreover, in his view, this manner was intended to disguise his own algebraic method described in the article of 13 April 1702[29]. He also pointed to other supposed technical errors that we do not discuss here. His remarks showed the depth to which he did not accept that Transcendental Geometry should adopt simple rules of algebra to treat differentials as if they were magnitudes, e.g. the rule of simply substituting equal things for equal things.

Saurin's reply of 23 April 1705 reiterated each of Rolle's criticisms in order to emphasise their unfoundedness[30]. But, since the beginning of the conflict, the authorities of the Academy had not positioned themselves with regard to the dispute. In a letter of 10 May 1705[31], Varignon reminded Leibniz that in spite of "tout ce que je luy demontrois [à Rolle] de paralogismes commis par sa seule ignorance de votre calcul", he had not been able to obtain justice from the Commissioners nominated on 3 September 1701[32]. He accused Rolle of misleading the public by his "hardiesse" and Gallois by the "bruit qu'il fait en sa faveur", so that the judges – Bignon, Cassini, la Hire and Gouye – did not have the courage to decide. In these conditions, he claimed that it would only be possible to counteract Rolle by attestations from scholars. He and Saurin therefore requested him to send them not only his attestation but also that of other scholars who were connoisseurs of differential calculus. Varignon asked him to abstain from "parler d'infinis" and to focus on the points at the end of the article.

By pointing out not to "parler d'infinis", Varignon was referring to the discussions that had taken place on the subject of publications in the *Journal de Trévoux* four years earlier. In May 1701, Thomas Gouye (1650-1725) published an article critical of the differential calculus in which he attacked in particular the consideration of infinities of any order[33]. In a reply article, Leibniz explained that the practice of his calculus did not require taking the infinitesimally small in a "rigorous" way, but only as small quantities as necessary to make the error less than wished[34]. He argued that in computational practice it is possible to consider an infinitesimal to be relative to a finite quantity such as a ball to the diameter of the Earth or the Earth to the distance of the fixed points. Similarly, an infinitely infinitely small quantity is like the ball to the distance to the fixed points. These comparisons are very similar to those he used while explaining his famous Lemma of Incomparables[35].This "lemma" assumes that it is possible to neglect some quantities over others if the mathematician calculating with them judges that they should be considered "incomparable". For Leibniz, the incomparable quantity is that which can be considered as small as to make the error smaller than any given error[36].

---

[28] Rolle: *Remarques* (cf. footnote 5), p. 131.
[29] Rolle: "Règles et remarques" (cf. footnote 19), pp. 240–242.
[30] Saurin: "Défense de la réponse à M. Rolle" (cf. footnote 6).
[31] Letter of Pierre Varignon to Leibniz, 10 May [1705], A III, 9, 549 (Vorausedition, 2.11.2018).
[32] Cf. footnote 2.
[33] Thomas Gouye: « Nouvelle methode pour déterminer aisément les rayons de la développée dans toute sorte de courbe algébraïque. Par Monsieur Jacques B. Acta Eruditorum, Mensis Novembris anni 1700. Lipsiae », in: *Journal de Trévoux*, mai – juin 1701, pp. 422–430.
[34] Leibniz: « Mémoire de Mr Leibnitz » (cf. footnote 3).

The academicians – whether attacking or defending the differential calculus – were quite astonished that Leibniz should claim that a differential is a fixed and determinate quantity, because this was how they interpreted the comparisons with finite quantities introduced in the *Journal de Trévoux* article. Leibniz reformulated his arguments in a reply published in March 1702[37] which completely calmed the doubts of the academicians concerned[38]. Faced with these antecedents, Varignon above all wanted to prevent another episode of this type, as his imperative warning "sans parler d'infinis" shows.

Thus, in the memoir that accompanied the letter, Varignon took care to formulate six very precise demands. They concerned the strictly technical criticisms that Rolle had developed in articles XI and XIV of *Remarques*, indicating precisely the pages of the article where these refutations occurred. He added that it would also be appropriate to announce "quelque chose dans les actes de Leipsik".

Leibniz promptly reacted to Varignon's request as soon as he received Saurin's article of 23 April 1705. In the letter of 27 July[39], he explained that although he "n'aime pas les contestations", he had written an attestation. However, he also decided on his own to write to Bignon and Gallois[40]. He would like the attestation to be transmitted to Abbé Bignon, but not to be published at first, and then, if nothing was obtained in this way, to renounce academic recognition and solicit the République des lettres.

4. *Puisque des personnes que j'estime beaucoup*

Leibniz sent two letters at the end of July, one to Bignon, the other to Gallois. The letter to Bignon is shorter and Leibniz did not elaborate at length on mathematical or metamathematical arguments. There are, however, similarities. For example, Leibniz discussed the criticism of his calculus for its lack of apodictic value. This comment is particularly directed at Gallois, whom Leibniz described in his letter as "zélé" for the "conservation de la Geometrie dans son ancienne pureté". Leibniz pointed out that this is indeed the lot of new methods, which by their novelty are susceptible to justification. In the letter to Bignon, to illustrate this, he referred to Cavalieri's Geometry of Indivisibles and then to Descartes' analysis, which, like his calculus, was criticised for not being "assez demonstrative". By citing these examples of prestige, which history shows to have been finally recognised for their apodictic value, Leibniz probably wanted to suggest that this would also be the case for his calculus. Moreover, he himself was on the side of rigour, admiring mathematicians who, like Proclus or Roberval, "wanted to demonstrate Euclid's axioms[41]". He considered it useful to demonstrate in the manner of the Ancients all the discoveries that his calculus had produced. In any case, this kind of criticism could only help, according to Leibniz's optimistic nature, since it gave "occasion to some clarification".

---

[35] Leibniz: "Tentamen de motuum coelestium causis", in: *Acta Eruditorum*, February 1689, pp. 82–96, GM VI, 144–161.
[36] "Itaque si quis nolit adhibere *infinite parvas,* potest assumere tam parvas quam sufficere judicat, ut sint incomparabiles, & errorem nullius momenti, imo dato minorem, producant", Leibniz: "Tentamen" (cf. footnote 35), p. 85, GM VI, 151.
[37] Leibniz: "Extrait d'une lettre de M. Leibnitz à M. Varignon" (cf. footnote 3).
[38] *Der Briefwechsel von Johann I. Bernoulli*, im Auftr. d. Naturforschenden Gesellschaft in Basel u. d. Otto-Spiess-Stiftung, 1–(3) vols., Basel 1955–1992, Bd. 1: *Der Briefwechsel mit Jakob Bernoulli, dem Marquis de l'Hôpital und anderen*, hrsg. von David Speiser, Basel 1955, p. 235. Dorénavant « *DBJB*, 1 ». The next tomes concern the correspondence with Pierre Varignon: *Der Briefwechsel Johann I. Bernoulli*, Bd. 2 : *Der Briefwechsel mit Pierre Varignon. Erster Teil : 1692–1702*, bearb. und kommentiert von Pierre Costabel ; Jeanne Peiffer, unter Benutzung von Vorarb. von Joachim Otto Fleckenstein, Basel 1988, et *Der Briefwechsel Johann I. Bernoulli*, Bd. 3 : *Der Briefwechsel mit Pierre Varignon. Zweiter Teil : 1702–1714*, bearb. und kommentiert von Pierre Costabel und Jeanne Peiffer, Basel 1992. We refer as « *DBJB*, 2 » et « *DBJB*, 3 » followed by the page. Here: *DBJB*, 2, 312.
[39] Letter from Leibniz to Pierre Varignon, 27 July 1705, A III, 10, 32 (Vorausedition, 7.10.2020), GM IV, 127–128.
[40] *Ibid*. The letters to Bignon and Gallois are respectively consultable in A I, 24, 837, and A III, 10, 46–55 (Vorausedition, 7.10.2020).

However, in the case of Rolle's accusations against Saurin, Leibniz was peremptory: the criticisms addressed to Saurin were absolutely not of this kind. They fall under the ordinary rules of mathematics that the algebraic calculus has partly borrowed from the geometrical axioms and from which the differential calculus, as a calculus, inherits[42].

Thus, Leibniz asked Bignon and Gallois to condemn these unfortunate objections that opposed to evidential truths, insisting, especially to Bignon, on the damaging consequences that a non-judicial intervention would have on the harmony of the institutional community and the public in general:

> [I]l s'agit seulement de savoir si on ne l'a pas combattu icy d'une maniere tout à fait indigne. Et s'il n'est de l'interest de la justice et des sciences, qu'un si mauvais exemple ne passe pas, apres l'establissement d'une societé si autorisée, et sous ses yeux.
> 
> Le public vous aura l'obligation d'avoir arresté un abus, capable d'infecter les sciences mathematiques qui en ont esté les plus exemtes[43].

Although Leibniz advanced an argument in order to convince the President of the Academy, it was not, in my opinion, a mere strategy. Leibniz was not necessarily against scientific disputes, on the contrary, he considered that they often led to discussions beneficial for the advancement of science. However, he considered the dispute with Rolle was not just illegitimate, but abusive[44].

The letter to Gallois is especially interesting. Leibniz had corresponded with Gallois between 1672 and 1692 and this letter breaks a long interruption. Leibniz returned to the difficulties raised by his reply to the *Journal de Trévoux* in November-December 1701. As we have seen, this reply had created confusion among the members of the Academy, and it was for this reason that Varignon had expressly asked him not to use the ideas he had developed there. But Leibniz insisted on making geometrical sense of differentials by introducing them as "modifications" of finite quantities. Since these modifications proceed in an unassignable way, continuously and without making any leap, they can be continuously unequal and thus leading to increases called second differences, and so on, from "degree" to "degree". He emphasized that these "degrees of magnitudes" are heterogeneous to each other but nevertheless determine finite magnitudes. Thus, the differences of the ordinates determine the direction of a curve, the second differences determine the angle of contact or curvature, and the following ones determine the angles of osculation which are of all degrees. Furthermore, these degrees of magnitude are incomparable to each other but have the same dimension. Thus, the elimination of incomparable quantities receives a justification which is not purely logical, but based on a classification of geometric problems by means of the different degrees of the quantities involved.

---

[41] Letter from Leibniz to Jean Gallois, July1705 [end], A III, 10, 47 (Vorausedition, 7.10.2020).
[42] "Elles reviennent à dire en effect qu'en maniant ce nouveau Calcul des infinitesimales, on ne doit point avoir la liberté d'y joindre les axiomes et opérations de la Geometrie et de l'Analyse ancienne, qu'on ne doit pas substituer *aequalibus aequalia*, qu'on ne doit point dire que deux quantités égales les quarrés sont egaux aussi, et des choses semblables"; A I, 24, 838. Similar remarks to Gallois", A III, 10, 48 (Vorausedition, 7.10.2020).
[43] A I, 24, 839–840.
[44] In the letter of 4 September 1699, he explained to James Cressett: « Vous avés bien raison de dire que les contestations entre les personnes studieuses sont utiles, et servent à l'avancement des sciences : mais quand on n'y observe pas les loix de l'honnesteté, elles nuisent plus qu'elles ne servent, et decrient les sciences, en faisant voir, que ceux qui s'y adonnent n'en deviennenet point meilleurs », A I, 17, 439. Concerning Leibniz and controversies, see Marcelo Dascal (ed.): *The Practice of Reason. Leibniz and His Controversies*, Philadelphia 2010.

By making distinctions involving the notions of dimension and heterogeneity, Leibniz provided new clarifications to the lemma of incomparables that he so often invoked to justify the validity of his calculus. This is particularly interesting since, in the same years, while addressing the Parisian public, Leibniz had abandoned the argument according to which the rules of his calculus, based on the lemma of incomparables, could be justified in the manner of the Ancients. From 1702 onwards, he substituted justifications based on the principle of continuity[45], which had the advantage of making fewer considerations about the status of differentials.

As demanded by Varignon, Leibniz wrote an attestation to Bignon. The French and Latin versions of the attestation are almost identical; it was the Latin version that was published. Both are structured in two distinct parts. In the first part, Leibniz settled some technical aspects of Rolle's criticism. In the second part, he exposed the conceptual scope of article 163 which makes the text interesting from an epistemological point of view[46].

For ease of reading, Leibniz noted by $v$ $3y^2 - 12y - 2x + 16$, and by $z$, $y^3 - 6yy + 8y - 6xy + 12x$, so that $\frac{dy}{dx} = \frac{3y^2-12y-2x+16}{y^3-6yy+8y-6xy+12x} = \frac{v}{z}$.

Leibniz aimed to justify why $v$ and $z$ should not be interpreted as "nothing" but as "evanescent" when $x$ and $y$ become equal to 2, *i.e.* when this quotient becomes the expression "$\frac{0}{0}$". Leibniz wanted to show the superiority of his calculus which, unlike ordinary algebra, is able to interpret and process this expression. Indeed, Leibniz argued $\frac{0}{0}$ is indeterminable if one "remains in these terms", that is to say if one remains within the framework of algebra in which 0 has its characteristic properties, in particular for multiplication. He explained that each of the equalities

$$\frac{0}{0} = 1 \; ; \; \frac{0}{0} = 2 \; ; \ldots \frac{0}{0} = b$$

has no more reason to be true than any of the others, because multiplying each member of the equation by 0, they become

$$0 = 1 \times 0 \, , \, 0 = 2 \times 0, \ldots 0 = b \times 0$$

which are all algebraically true because of the algebraic nature of 0: "whatever quantity or number can be multiplied by 0, it always comes from 0". Thus, Leibniz concluded that it was not possible to decide whether the quotient is equal to 1, 2 or more generally $b$. By taking $v$ and $z$ to be nothing, "we learn nothing" [*ex quibus nihil discemus*] from this relationship. The algebraic framework therefore does not allow us to give meaning to the expression $\frac{0}{0}$, a meaning that is nevertheless indispensable to solve the problem of the double point of the quartic curve.

---

[45] Rabouin/Arthur: "Leibniz's Syncategorematic Infinitesimals II" (cf. footnote 4).
[46] This text has been analysed in Sandra Bella: "De l'intraitable à l'indéterminé: entre calcul et géométrie, réflexions leibniziennes autour de $\frac{0}{0}$ (1700–1705)", in: *Mathématique et philosophie leibniziennes à la lumière des manuscrits inédits* (*Cahier thématique de Philosophia Scientiae* 25/2), éd. par Valérie Debuiche et David Rabouin, June 2021. The main interests are summarised here.

It is noteworthy that Leibniz did not state here that $v$ and $z$ are not nothing, but that if we want to perform a computing operation, we must not take them for nothing. To have to consider them otherwise than as nothings is thus an imperative linked to a calculatory practice: "they [$v$ and $z$] must be taken as nascent or vanishing magnitudes". Incipient or vanishing magnitudes are those which differ from nothing in an unassignable way. Although they are not assignable, they are something existing and not an absolute nothing[47]. Moreover, the magnitude $z$, incipient or vanishing, is "as much as its $dz$" because $z$ differing from nothing only in an unassignable way[48], its first or ultimate infinitesimal difference, existing at the same moment, coincides with it[49]. This last comment justifies the final result of Article 163, which orders that instead of the quotient of the numerator and the denominator, the quotient of their differentials should be taken.

These considerations allowed Leibniz to conclude that only differential calculus, through the notion of differentials/evanescents, could provide an interpretation of the expression $\frac{0}{0}$, which had occurred in the double point configuration. The explanation of article 163 thus shows him the superiority of his calculus over the ordinary algebra.

5. Publications and outcomes

According to Jacques Lelong[50] and Varignon[51], Leibniz's written intervention influenced Bignon who on 8 August 1705 named a commission to judge the conflict between Rolle and Saurin. The commission was composed of Gallois (director), Cassini, Gouye, Fontenelle and La Hire[52]. Gouye and La Hire were not very favourable to the introduction of the new calculus, not to speak of Gallois. Such a choice of composition could hardly be expected to result in a verdict favourable to the camp of the calculus defenders. Shortly after their appointment, each of them sent Bignon his verdict in writing, but the latter had not yet made it public, because, according to Lelong, "il parroist qu'il vaut sauver l'honneur de M$^r$ Rolle[53]".

Finally, on 9 January 1706, Fontenelle read the Commission's judgment[54]. According to the procès-verbal, the verdict faithfully summarised all the sentiments expressed secretly by the five members in fifteen articles. First of all, Rolle was exonerated: when he stated that the "ordinary methods" do not allow to solve the general problem of tangents, he did not necessarily include the differential calculus among these methods. His aim was only to show that his own algebraic method was able to solve the problem of tangents in a general way.

---

[47] "ubi adhuc aliquid aut jam aliquid rei adest, etsi assignabile non sit", Gotha FB A 448–449, fol. 42r◦ ; Leibniz: « Sentiment de Monsieur Leibnitz » (cf. note 3), p. 3.
[48] "[…] a nihilo inassignabiliter differentibus", ibid.
[49] "Porro quantitas $z$, ubi est nascens aut evanescens, idem est quod $dz$ ibidem existens ; […] cum magnitudo nascens aut evanescens, seipsa a nihilo differat inassignabiliter, primaque adeo seu ultima differentia (infinitesimalis scilicet) coincidat cum termino ipso.", ibid.
[50] Letter from Lelong to Leibniz, 5 October 1705, A I, 25, 110.
[51] Letter from Varignon to Leibniz, 9 October 1705, A III, 10, 100 (Vorausedition, 7.10.2020), GM IV, 131.
[52] *PVARS*, 24, 263v◦-264r◦.
[53] A I, 25, 110.
[54] *PVARS*, 25, 1r◦-4r◦.

Thus, to contradict Rolle's statement, Saurin should have shown the generality of his method, derived from the differential calculus. In the last article, the Commission asked the two opponents to adopt a conduct worthy of scholars and honest men[55].

The decision was highly awaited, especially by the defenders who were expecting to obtain compensation. But, the outcome of this episode was motivated more by the wish for institutional appeasement than by granting a victory to differential calculus.

In the meantime, three other attestations were included in the 1706 publication. Those of Hermann and Jacob Bernoulli were written on May (respectively 21 and 26 May) in Latin. Hermann's and Jacob Bernoulli's attestations only ratified the six points of Varignon's memoir in order to reject Rolle's attacks.

Johann Bernoulli's attestation, written in French, did not arrive until the first trimester of 1706, whereas Varignon had already requested it on 30 May[56]. This delay can be explained[57]. Although Johann Bernoulli structured his writing by dealing one by one with the six points of Varignon's memoir, sarcastic comments were added about Rolle's intellectual blindness[58].

Saurin published his article "CONTINUATION DE LA DEFENSE de M. Saurin contre la Replique de M. Rolle publié en 1703, sous le titre de Remarques touchant le Problême général des Tangentes, &c" with the four attestations he entitled "Sentimens de quelques Geometres celebres, sur le premier article de cette dispute".

The whole was printed in "Amsterdam, chés Henry Westein" but several indices would tend to show that the printing was probably realized clandestinely in France[59].

In a letter to Leibniz dated 8 March 1706, Lelong, without detailing all the points of the academic judgement, mentions Saurin's publication and his peculiar destiny. Indeed, on discovering that Saurin had taken the liberty of such an audacious act, Bignon confiscated all the copies and had threatened to revoke his pension[60]. Moreover, no French newspaper reported on the verdict. The January issue of the *Nouvelles de la République des Lettres* – which Leibniz read[61]– provided a very brief summary of the verdict and in fact only included the last article of the academic judgment[62]. Thus, Leibniz was largely unaware of the details of the judgement, but despite this lack of information he surmised that Bignon's aim was appeasement rather than a true verdict on mathematical results: "magis morale quam mathematicum", he concluded to Johann Bernoulli[63].

---

[55] "Mais si dans cette pensée, ou par telle raison que ce puisse être, ces Mrs continüet leur dispute, ils ne sauroient se repñter trop souvent et trop fortement, Les regles que toutes sortes de Loix imposet. à des Savans cõme à d'honêtes gens ; Le public n'ayant été que trop scandalisé de voir que dans leurs divers Ecrits ils se font abandonez à trop d'aigreur & à des manières trop méprisantes de part et d'autre, ce que souvent même ils y ont trop dôné à l'esprit aux dépens de la politesse, de la droiture, et de la bône foy. C'est surquoy Nous renvoyons Mr Rolle aux Statuts de l'Académie, dont il a l'honeur d'être, & Mr Saurin à son propre coeur", PVARS, 25, 4r°.
[56] *DBJB*, 3, 151–152.
[57] Bella: *La Reconstruction* (cf. footnote 19), p. 450-452
[58] " […] qu'il faut être aussi aveugle que M. Rolle pour ne s'en apercevoir pas d'abord", Johann Bernoulli: "Déclaration de J. Bernoully sur les Articles proposés par M. Saurin, concernant la dispute qui est entre lui & M. Rolle", in: Saurin: *Continuation de la défense* (cf. footnote 3), pp. 40–42, p. 40, "Non seulement c'est une erreur grossiere, mais aussi une opiniâtreté insuportable ; M. Rolle attaque des choses qu'il n'entend pas, c'est un aveugle qui parle des couleurs », ibid., p. 41, « mais elle ne sera manifeste aux aveugles & aux entêtés tels que M. Rolle", ibid., p. 41.
[59] Malebranche had also published a work under the name of the printer 'Henry Westein' in 1704 [Réponse du Père Malebranche, à la troisième lettre de M. Arnauld, 1704], but the printing was actually made in Rouen. In fact, there is a printer in Amsterdam with the name "Wetstein" (*i.e.* Westein with an extra letter). It is therefore very likely that "Westein" is simply a figurehead for a printer in Rouen. It is also possible, but not certain, that the printer of Malebranche's work is the same as that of Saurin. We thank Dr. Charlotte Wahl and Dr. Siegmund Probst for these valuable indications.
[60] A I, 25, 702.
[61] GM III, 789.
[62] "les Commissaires nommez pour examiner le différent qui étoit entre Mr Saurin et Mr Rolle sur les infiniment petits, prononcèrent leur jugement, & renvoyèrent dit-on Mr Rolle, aux statuts de l'Académie qui ordonnent qu'on dira les choses avec ménagement ; & à l'égard de Mr Saurin à son bon cœur", *Nouvelles de la République des Lettres*, janvier 1706, p. 120.

On 19 August 1706, Leibniz wrote directly to Bignon. He informed him that he had learned of the Academy's "sort of decision" from the Dutch newspaper. Without blaming Bignon, Leibniz explained that "good emulation" consisted in either re-demonstrating results in the manner of the ancients or inventing new methods or improving the present ones[64]. Although Leibniz did not write that this was exactly what Rolle did not do, the implication was patent. Bignon's reply was as sincere as his role as President of the Academy would have allowed:

> "Je suis persuadé comme vous que la dispute qui etoit entre M$^{rs}$ Saurin et Rolle, ne pouvoit aboutir à rien davantageux et l'academie a pris un fort bon parti en la faisant cesser sans se commettre.[65]"

This tempered position probably made Leibniz give up insisting on his point of view to Bignon and made him conclude in the next letter:

> "Il est bon au moins que la dispute entre vos Geomestres a esté terminée sans que la verité en ait souffert[66]."

6. Conclusion

Leibniz was the inventor of a completely new method of calculus which he made public. This publicity implied proving its validity to his peers. This was not an easy task, as he had to convince them that his method followed a continuum of practice and that it was equivalent to an old-style demonstration. Probably because of his diplomatic qualities, Leibniz was sensitive to everyone's affinities when he interceded to defend his calculus. This is evidenced by his public interventions as well as his epistolary exchanges, the two letters of which analysed here are particularly significant.

When Varignon asked him to mediate, he wanted a well-ordered intervention, which he oriented. However, Leibniz did not follow this prescription. His attestation constitutes a justification of his calculus that would be comparable to his contemporary texts (*Justification du Calcul des infinitésimales par celuy de l'Algèbre ordinaire et Defense du calcul des Differences*). The differences between the *analyse des infinis* and ordinary analysis is that the former introduces entities that allow computations to be made, whereas the latter inexorably fails. In "Sentiment de Monsieur Leibnitz", these entities are presented as "evanescent magnitudes", a terminology not unlike Newton's Principia[67]. It is not a matter of being concerned with the existence of these entities, but of "taking" them for something other than 'pure nothings'. The strategy of appealing to the notion of evanescents is the one he definitely chose with the Parisian public. But in the letter to Gallois, in a private way, his justification took another form which appealed to the notion of the incomparable.

It is quite clear that the different notions used in his justifications allowed Leibniz to adapt his discourse. But it seems to me that what appears to be hesitations should rather be interpreted as a pursuit –certainly not completed during his lifetime– of a coherent theory for the foundation of his calculus, beyond his interlocutors, in which the various notions invoked – the incomparable, the unassignable, the evanescent – can each be understood from the others.

---

[63] GM III, 794.
[64] Letter from Leibniz to Jean Paul Bignon, 19 August 1706, A I, 26, 406 (Vorausedition, 25.9.2020).
[65] Letter from Jean Paul Bignon to Leibniz, 13 October 1706, A I, 26, 574 (Vorausedition, 25.9.2020).
[66] Letter from Leibniz to Jean Paul Bignon, 23 April 1707, A I, 27, 182 (Vorausedition, 13.11.2020).
[67] Arthur: « Leery Bedfellows » (cf. footnote 4), p. 19.

The analysis of this episode illustrates that, in general, the acceptance of the Leibnizian calculus in the Académie royale des sciences did not depend exclusively on logical arguments or on criteria internal to the discipline, that it was not necessary to share a philosophy of mathematics either, but that the achievement of a consensus was the result of philosophical but also strongly political debates.